\documentclass[12pt,a4paper]{article}
\usepackage{amsthm,amsmath,amssymb}
\theoremstyle{plain}
\newtheorem{Theorem}{Theorem}[section]
\newtheorem{Lemma}[Theorem]{Lemma}
\newtheorem{Proposition}[Theorem]{Proposition}
\theoremstyle{definition}
\newtheorem{Definition}[Theorem]{Definition}
\newtheorem{Remark}[Theorem]{Remark}
\newcommand{\NNN}{\mathbb N}
\newcommand{\gen}[1]{\mathopen{<}#1\mathclose{>}}
\newcommand{\abs}[1]{\lvert#1\rvert}
\DeclareMathOperator{\GL}{GL}
\DeclareMathOperator{\PGL}{PGL}
\DeclareMathOperator{\Gal}{Gal}
\begin{document}
\title{Quadratic Factors of $f(X)-g(Y)$}
\author{Manisha Kulkarni \and Peter M\"uller \and B.~Sury}
\maketitle
\section{Introduction}
This note extends the characteristic $0$ results in \cite{Bilu:qf} to
arbitrary characteristic. The method is completely different from Bilu's.  The
main bulk of the work handles the case of positive characteristic. Indeed, if
one skips all the arguments specific to this, one obtains a particularly short
and natural proof of Bilu's results. Also, the rather specific main result of
\cite{BerrondoGallardo} is a trivial consequence of the theorems below.

The generalization of \cite[Theorem 1.2]{Bilu:qf} is
\begin{Theorem}\label{T:main1}
  Let $f,g\in K[X]$ be polynomials over a field $K$, such that $f(X)-g(Y)\in
  K[X,Y]$ has a factor of degree at most $2$. If the characteristic
  $p$ of $K$ is positive, then assume that $f$ or $g$ cannot be written as a
  polynomial in $X^p$. Then there are $f_1,g_1,\Phi\in K[X]$ with $f=\Phi\circ
  f_1$, $g=\Phi\circ g_1$, such that one of the following holds:
  \begin{itemize}
  \item[(a)] $\deg f_1,\deg g_1\le 2$.
  \item[(b)] $p\ne2$, $n=\deg f_1=\deg g_1\ge 4$ is a power of $2$, and there
    are $\alpha,\beta,\gamma,a\in K$ such that $f_1(X)=D_n(X+\beta,a)$,
    $g_1(X)=-D_n(\alpha X+\gamma(\xi+1/\xi),a)$. Here $\xi$ denotes a
    primitive $2n$-th root of unity. Furthermore, $\xi^2+1/\xi^2\in K$.
  \end{itemize}
\end{Theorem}
Conversely, in cases (a) and (b) $f(X)-g(Y)$ indeed has a factor of degree at
most $2$. This is clear for case (a), because $f_1(X)-g_1(Y)$ is such a
factor, and follows for case (b) from Lemma \ref{L:dicksonfact}.

If one wants to determine the cases such that $f(X)-g(Y)$ has an irreducible
factor of degree $2$, then the list becomes longer in positive characteristic.
The exact extension of \cite[Theorem 1.3]{Bilu:qf} is
\begin{Theorem}\label{T:main2}
  Let $f,g\in K[X]$ be polynomials over a field $K$, such that $f(X)-g(Y)\in
  K[X,Y]$ has a quadratic irreducible factor $q(X,Y)$. If the characteristic
  $p$ of $K$ is positive, then assume that $f$ or $g$ cannot be written as a
  polynomial in $X^p$. Then there are $f_1,g_1,\Phi\in K[X]$ with $f=\Phi\circ
  f_1$, $g=\Phi\circ g_1$ such that $q(X,Y)$ divides $f_1(X)-g_1(Y)$, and one
  of the following holds:
  \begin{itemize}
  \item[(a)] $\max(\deg f_1,\deg g_1)=2$ and $q(X,Y)=f_1(X)-g_1(Y)$.
  \item[(b)] There are $\alpha,\beta,\gamma,\delta\in K$ with
    $g_1(X)=f_1(\alpha X+\beta)$, and $f_1(X)=h(\gamma X+\delta)$, where
    $h(X)$ is one of the following polynomials.
  \begin{itemize}
  \item[(i)] $p$ does not divide $n$, and $h(X)=D_n(X,a)$ for some $a\in K$.
    If $a\ne0$, then $\zeta+1/\zeta\in K$ where $\zeta$ is a primitive $n$-th
    root of unity.
\item[(ii)] $p\ge3$, and $h(X)=X^p-aX$ for some $a\in K$.
\item[(iii)] $p\ge3$, and $h(X)=(X^p+aX+b)^2$ for some $a,b\in K$.
\item[(iv)] $p\ge3$, and $h(X)=X^p-2aX^{\frac{p+1}{2}}+a^2X$ for some $a\in K$.
\item[(v)] $p=2$, and $h(X)=X^4+(1+a)X^2+aX$ for some $a\in K$.
\end{itemize}
\item[(c)] $n$ is even, $p$ does not divide $n$, and there are
  $\alpha,\beta,\gamma,a\in K$ such that $f_1(X)=D_n(X+\beta,a)$,
  $g_1(X)=-D_n(\alpha X+\gamma(\xi+1/\xi),a)$. Here $\xi$ denotes a primitive
  $2n$-th root of unity. Furthermore, $\xi^2+1/\xi^2\in K$.
\item[(d)] $p\ge 3$, and there are quadratic polynomials $u(X),v(X)\in K[X]$,
  such that $f_1(X)=h(u(X))$ and $g_1(X)=h(v(X))$ with
  $h(X)=X^p-2aX^{\frac{p+1}{2}}+a^2X$ for some $a\in K$.
\end{itemize}
\end{Theorem}
The theorems exclude the case that $f$ and $g$ are both polynomials in $X^p$.
The following handles this case, a repeated application reduces to the
situation of the Theorems above.
\begin{Theorem}\label{T:X^p}
  Let $f,g\in K[X]$ be polynomials over a field $K$, such that $f(X)-g(Y)\in
  K[X,Y]$ has an irreducible factor $q(X,Y)$ of degree at most $2$. Suppose
  that $f(X)=f_0(X^p)$ and $g(X)=g_0(X^p)$, where $p>0$ is the characteristic
  of $K$. Then one of the following holds:
  \begin{itemize}
  \item[(a)] $q(X,Y)$ divides $f_0(X)-g_0(Y)$, or
\item[(b)] $p=2$, $f(X)=f_0(X^2)$, $g(X)=f_0(aX^2+b)$ for some $a,b\in K$, and
  $q(X,Y)=X^2-aY^2-b$.
  \end{itemize}
\end{Theorem}
\begin{Remark}
  Under suitable conditions on the parameters and the field $K$, all cases
  listed in Theorem \ref{T:main2} give examples such that $f_1(X)-g_1(Y)$
  indeed has an irreducible quadratic factor. The cases of the Dickson
  polynomials are classically known, see Lemma \ref{L:dicksonfact} and its
  proof. We illustrate two examples:

(b)(v). Here $p=2$ and $h(X)=X^4+(1+a)X^2+aX$. We have
$h(X)-h(Y)=(X+Y)(X+Y+1)(X^2+X+Y^2+Y+a)$. If $Z^2+Z=a$ has no solution in $K$,
then the quadratic factor is irreducible.

(b)(iv). Here $p\ge3$ and $h(X)=X^p-2aX^{\frac{p+1}{2}}+a^2X$, and $a\ne0$ of
course. If $\alpha$ is a root of $Z^{p-1}-a$, then so is $-\alpha$. Let $T$ be
a set such $T\cup(-T)$ is a disjoint union of the roots of $Z^{p-1}-a$.

We compute
\begin{align*}
  h(X^2)-h(Y^2) &= (X^2-Y^2)\prod_{t\in T\cup(-T)}[((X-Y)-t)((X+Y)-t)]\\
&=(X^2-Y^2)\prod_{t\in T}[((X-Y)-t)((X+Y)-t)\\
& \phantom{(X^2-Y^2)\prod_{t\in T}888}((X+Y)+t)((X-Y)+t)]\\
&=(X^2-Y^2)\prod_{t\in T}((X^2-Y^2)^2-2t^2(X^2+Y^2)+t^4).
\end{align*}
and therefore
\[
h(X)-h(Y)=(X-Y)\prod_{t\in T}((X-Y)^2-2t^2(X+Y)+t^4).
\]
The discriminant with respect to $X$ of the quadratic factor belonging to $t$
is $16t^2Y$, so all the quadratic factors are absolutely irreducible.
\end{Remark}
\section{Preparation}
\begin{Definition}
  Let $a,b$ elements of a group $G$. Then $a^b$ denotes the conjugate
  $b^{-1}ab$.
\end{Definition}
\begin{Lemma}\label{L:dieder2}
  Let $G$ be a finite dihedral group, generated by the involutions $a$ and
  $b$. Then $a$ and a suitable conjugate of $b$ generate a Sylow $2$-subgroup
  of $G$.
\end{Lemma}
\begin{proof}
  Set $c=ab$. For $i\in\NNN$, the order of $\gen{a,b^{c^i}}$ is twice the
  order of $ab^{c^i}$. We compute
  $ab^{c^i}=a(c^{-1})^ibc^i=a(ba)^ib(ab)^i=(ab)^{2i+1}=c^{2i+1}$. Let $2i+1$
  be the largest odd divisor of $\abs{G}$. The claim follows.
\end{proof}
\begin{Definition}
  For $a,b,c,d$ in a field $K$ with $ad-bc\ne0$ let
  $\begin{bmatrix}a&b\\c&d\end{bmatrix}$ denote the image of
  $\begin{pmatrix}a&b\\c&d\end{pmatrix}\in\GL_2(K)$ in $\PGL_2(K)$.
\end{Definition}
\begin{Lemma}\label{L:pglcyc}
  Let $K$ be an algebraically closed field of characteristic $p$, and
  $\rho\in\PGL_2(K)$ be an element of finite order $n$. Then one of the
  following holds:
  \begin{itemize}
  \item[(a)] $p$ does not divide $n$, and $\rho$ is conjugate to
    $\begin{bmatrix}1&0\\0&\zeta\end{bmatrix}$, where $\zeta$ is a primitive
    $n$-th root of unity.
  \item[(b)] $n=p$, and $\rho$ is conjugate to
    $\begin{bmatrix}1&1\\0&1\end{bmatrix}$.
  \end{itemize}
\end{Lemma}
\begin{proof}
  Let $\hat\rho\in\GL_2(K)$ be a preimage of $\rho$. Without loss of
  generality we may assume that $1$ is an eigenvalue of $\hat\rho$. The claim
  follows from the Jordan normal form of $\hat\rho$.
\end{proof}
\begin{Lemma}\label{L:pgldih}
  Let $K$ be an algebraically closed field of characteristic $p$, and
  $G\le\PGL_2(K)$ be a dihedral group of order $2n\ge4$, which is generated by
  the involution $\tau$ and the element $\rho$ of order $n$. Then one of the
  following holds:
  \begin{itemize}
  \item[(a)] $p$ does not divide $n$. There is $\sigma\in\PGL(K)$ such that
    $\tau^\sigma=\begin{bmatrix}0&1\\1&0\end{bmatrix}$ and
    $\rho^\sigma=\begin{bmatrix}1&0\\0&\zeta\end{bmatrix}$, where $\zeta$ is a
    primitive $n$-th root of unity.
  \item[(b)] $n=p\ge3$. There is $\sigma\in\PGL(K)$ such that
    $\tau^\sigma=\begin{bmatrix}1&0\\0&-1\end{bmatrix}$ and
    $\rho^\sigma=\begin{bmatrix}1&1\\0&1\end{bmatrix}$.
  \item[(c)] $n=p=2$. There is $\sigma\in\PGL(K)$ such that
    $\tau^\sigma=\begin{bmatrix}1&b\\0&1\end{bmatrix}$ and
    $\rho^\sigma=\begin{bmatrix}1&1\\0&1\end{bmatrix}$ for some $1\ne b\in K$.
  \end{itemize}
\end{Lemma}
\begin{proof}
By Lemma \ref{L:pglcyc} we may assume that $\rho$ has the form given
there. From $\rho^\tau=\rho^{-1}$ we obtain the shape of $\tau$:

First assume that $p$ does not divide $n$, so
$\rho=\begin{bmatrix}1&0\\0&\zeta\end{bmatrix}$. Let
$\hat\tau=\begin{pmatrix}a&b\\c&d\end{pmatrix}\in\GL_2(K)$ be a preimage of
$\tau$. From $\rho^\tau=\rho^{-1}$ we obtain $\rho\tau=\tau\rho^{-1}$, hence
\[
 \begin{pmatrix}1&0\\0&\zeta\end{pmatrix}
\begin{pmatrix}a&b\\c&d\end{pmatrix}=
\lambda\begin{pmatrix}a&b\\c&d\end{pmatrix}
\begin{pmatrix}\zeta&0\\0&1\end{pmatrix}
\]
for some $\lambda\in K$. This gives $(\lambda\zeta-1)a=0$, $(\lambda-1)b=0$,
$(\lambda-1)c=0$, and $(\lambda-\zeta)d=0$. First assume $b=c=0$. Then $\rho$
and $\tau$ commute, so $G$ is abelian, hence $n=2\ne p$ and therefore
$\zeta=-1$. It follows $\tau=\begin{bmatrix}1&0\\0&-1\end{bmatrix}=\rho$, a
contradiction.

Thus $b\ne0$, so $\lambda=1$. This yields $a=d=0$, as $\zeta\ne1$. We obtain
$\tau=\begin{bmatrix}0&1\\c&0\end{bmatrix}$. Choose $\beta\in K$ with
$\beta^2=c$, and set $\delta=\begin{bmatrix}1&\beta\\0&1\end{bmatrix}$. The
claim follows from $\rho^\delta=\rho$ and
$\tau^\delta=\begin{bmatrix}0&1\\1&0\end{bmatrix}$.

Now assume the second case of Lemma \ref{L:pglcyc}, that is $p=n$ and
$\rho=\begin{bmatrix}1&1\\0&1\end{bmatrix}$. Again setting
$\hat\tau=\begin{pmatrix}a&b\\c&d\end{pmatrix}$ we obtain
\[
 \begin{pmatrix}1&1\\0&1\end{pmatrix}
\begin{pmatrix}a&b\\c&d\end{pmatrix}=
\lambda\begin{pmatrix}a&b\\c&d\end{pmatrix}
\begin{pmatrix}1&-1\\0&1\end{pmatrix}
\]
for some $\lambda\in K$. This gives $a+c=\lambda a$, $b+d=\lambda(-a+b)$,
$c=\lambda c$, and $d=\lambda(-c+d)$. If $c\ne0$, then $\lambda=1$, so $c=0$
by the first equation, a contradiction. Thus $c=0$, so $a\ne0$. We may assume
$a=1$, so $d=-1$. This gives the result for $p=n=2$. If $p\ne2$, then set
$\sigma=\begin{bmatrix}1&\beta\\0&1\end{bmatrix}$ with $\beta=-b/2$. From
$\rho^\sigma=\rho$ and $\tau^\sigma=\begin{bmatrix}1&0\\0&-1\end{bmatrix}$ we
obtain the claim.
\end{proof}
Let $z$ be a transcendental over the field $K$. The group of $K$-automorphisms
of $K(z)$ is isomorphic to $\PGL_2(K)$, where
$\begin{bmatrix}a&b\\c&d\end{bmatrix}$ sends $z$ to $\frac{az+b}{cz+d}$. Note
that $K(z)=K(z')$ for $z\in K(z)$ if and only if $z'=\frac{az+b}{cz+d}$ with
$\begin{bmatrix}a&b\\c&d\end{bmatrix}\in\PGL_2(K)$.

Let $r(z)\in K(z)$ be a rational function. Then the \emph{degree} $\deg r$ of
$r$ is the maximum of the degrees of the numerator and denominator of $r(z)$
as a reduced fraction. Note that $\deg r$ is also the degree of the field
extension $K(z)/K(r(z))$.
\begin{Definition}\label{D:dickson}
For $a\in K$ one defines the $n$th Dickson polynomial $D_n(X,a)$ (of degree
$n$) implicitly by $D_n(z+a/z,a)=z^n+(a/z)^n$. Note that $D_n(X,0)=X^n$.
Furthermore, from
$b^nD_n(z+a/z,a)=b^n(z^n+(a/z)^n)=(bz)^n+(\frac{b^2a}{bz})^n=
D_n(bz+\frac{b^2a}{bz},b^2a)=D_n(b(z+a/z),b^2a)$ one obtains
$b^nD_n(X,a)=D_n(bx,b^2a)$, a relation we will use later.
\end{Definition}
\begin{Lemma}\label{L:LRlr}
  \begin{itemize}
  \item[(a)] Let $f(X)=g(h(X))$ with $f\in K[X]$ and $g,h\in K(X)$. Then
    $f=g\circ h=(g\circ\lambda^{-1})\circ(\lambda\circ h)$ for a rational
    function $\lambda\in K(X)$ of degree $1$, such that $g\circ\lambda^{-1}$
    and $\lambda\circ h$ are polynomials.
  \item[(b)] Let $f,g\in K[X]$ be two polynomials such that $f(X)=L(g(R(X)))$
    for rational functions $L,R\in K(X)$ of degree $1$. Then there are linear
    polynomials $\ell,r\in K[X]$ with $f(X)=\ell(g(r(X)))$.
  \end{itemize}
\end{Lemma}
\begin{proof}
  (a) This is well known. For the convenience of the reader, we supply a short
  proof. Let $\lambda\in K(X)$ be of degree $1$ such that
  $\lambda(h(\infty))=\infty$. Setting $\bar g=g\circ\lambda^{-1}$ and $\bar
  h=\lambda\circ h$ we have $f=\bar g\circ\bar h$ with $\bar
  h(\infty)=\infty$. Suppose that $\bar g$ is not a polynomial. Then there is
  $\alpha\in\bar K$ ($\bar K$ denotes an algebraic closure of $K$) with $\bar
  g(\alpha)=\infty$. Let $\beta\in\bar K\cup\{\infty\}$ with $\bar
  h(\beta)=\alpha)$. From $\bar h(\infty)=\infty$ we obtain
  $\beta\ne\infty$. Now $f(\beta)=\bar g(\bar h(\beta))=\bar g(\alpha)=\infty$
  yields a contradiction, so $\bar g$ is a polynomial. From that it follows
  that $\bar h$ is a polynomial as well.

  (b) If $L$ is a polynomial, then $R$ has no poles, so is a polynomial as
  well.

  Suppose now that $L$ is not a polynomial. Then there is $\alpha\in K$ with
  $L(\alpha)=\infty$. Let $\bar K$ be an algebraic closure of $K$. Choose
  $\beta\in\bar K$ with $g(\beta)=\alpha$. If we can find $\gamma\in\bar K$
  with $R(\gamma)=\beta$, then we get the contradiction $f(\gamma)=\infty$.
  The value set of $R$ on $\bar K$ is $\bar K$ minus the element $R(\infty)\in
  K$. Thus we are done except for the case that the equation $g(X)=\alpha$ has
  only the single solution $\beta=R(\infty)\in K$. In this case, however,
  $g(X)=\alpha+\delta(X-\beta)^n$ with $\delta\in K$. From
  $L^{-1}(f(R^{-1}(X)))=g(X)$ we analogously either get that $L$ and $R$ are
  polynomials, or $f(X)=\alpha'+\delta'(X-\beta')^n$ with
  $\alpha',\delta',\beta'\in K$. The claim follows.
\end{proof}
\begin{Lemma}\label{L:dicksonfact}
  Let $K$ be a field of characteristic $p$, and $n\in\NNN$ even and not
  divisible by $p$ (so in particular $p\ne2$). Let $\xi$ be a primitive $2n$-th
  root of unity and $a\in K$. Then
\[
D_n(X,a)+D_n(Y,b)=\prod_{1\le k\le n-1\text{ odd}}
(X^2-(\xi^k+1/\xi^k)XY+Y^2-(\xi^k-1/\xi^k)^2a).
\]
\end{Lemma}
\begin{proof}
  This is essentially \cite[Prop.~3.1]{Bilu:qf}. The factorizations of
  $D_m(X,a)-D_m(Y,a)$ are known, see \cite[Prop.~1.7]{Turnwald:Schur}. The
  claim then follows from that and
  $D_{2n}(X,a)-D_{2n}(Y,b)=D_n(X,a)^2-D_n(Y,b)^2=
  (D_n(X,a)+D_n(Y,b))(D_n(X,a)-D_n(Y,b))$.
\end{proof}
The following proposition classifies polynomials $f$ over $K$ with a certain
Galois theoretic property. To facilitate the notation in the statement and its
proof, we introduce a notation: If $E$ is a field extension of $K$, and
$f,h\in K[X]$ are polynomials, then we write $f\sim_Eh$ if and only if there
are linear polynomials $L,R\in E[X]$ with $f(X)=L(h(R(X)))$. Clearly, $\sim_E$
is an equivalence relation on $K[X]$. In determining the possibilities of $f$
in Proposition \ref{P:poldih}, we first determine certain polynomials
$h\in\bar K[X]$ with $f\sim_{\bar K}h$, and from that we conclude the
possibilities for $f$. The following Lemma illustrates this latter step.
\begin{Lemma}\label{L:sim}
  Let $\bar K$ be an algebraic closure of the field $K$ of characteristic
  $p$. Suppose that $f\sim_{\bar K}X^p-2X^{(p+1)/2}+X$ for $f\in K[X]$. Then
  $f\sim_KX^p-2aX^{(p+1)/2}+a^2X$ for some $a\in K$.
\end{Lemma}
\begin{proof}
  There are $\alpha,\beta,\gamma,\delta\in\bar K$ with $f(X)=\alpha h(\gamma
  X+\delta)+\beta\in K[X]$, where $h(X)=X^p-2X^{(p+1)/2}+X$.

  The coefficients of $X^p$ and $X^{(p+1)/2}$ of $f(X)$ are $\alpha\gamma^p\in
  K$ and $-2\alpha\gamma^{(p+1)/2}\in K$, so $\gamma^{(p-1)/2}\in K$ and
  $\alpha\gamma\in K$.

  Suppose that $p>3$. Then the coefficient of $X^{(p-1)/2}$ is (up to a factor
  from $K$) $\alpha\gamma^{(p-1)/2}\delta\in K$, so $\alpha\delta\in K$ and
  therefore $\delta/\gamma\in K$. Thus, upon replacing $X$ by
  $X-\delta/\gamma$, we may assume $\delta=0$. Then $\beta\in K$, so $\beta=0$
  without loss of generality. Now dividing by $\alpha\gamma^p$ and setting
  $a=1/\gamma^{(p-1)/2}$ yields the claim.

  In the case $p=3$ we get from above $\gamma\in K$ and then $\alpha\in K$.
  Thus we may assume $\alpha=\gamma=1$. Looking at the coefficient of $X$,
  which is $-4\delta+1$, shows $\delta\in K$, so $\delta=\beta=0$ without loss
  of generality. Thus $f(X)=X^3-2X^2+X$.
\end{proof}
\begin{Proposition}\label{P:poldih}
  Let $K$ be a field of characteristic $p$, and $f(X)\in K[X]$ be a polynomial
  of degree $n\ge3$ which is not a polynomial in $X^p$. Let $x$ be a
  transcendental, and set $t=f(x)$. Suppose that the normal closure of
  $K(x)/K(t)$ has the form $K(x,y)$ where $F(x,y)=0$ with $F\in K[X,Y]$
  irreducible of total degree $2$. Furthermore, suppose that the Galois group
  of $K(x,y)/K(t)$ is dihedral of order $2n$. Then one of the following holds:
%  where $L,R\in K[X]$ are linear polynomials.
  \begin{itemize}
  \item[(a)] $p$ does not divide $n$, and $f\sim_KD_n(X,a)$ for some $a\in K$.
    If $a\ne0$, then $\zeta+1/\zeta\in K$ where $\zeta$ is a primitive $n$-th
    root of unity.
\item[(b)] $n=p\ge3$, and $f\sim_KX^p-aX$ for some $a\in K$.
\item[(c)] $n=2p\ge6$, and $f\sim_K(X^p+aX+b)^2$ for some $a,b\in K$.
\item[(d)] $n=p$, and $f\sim_KX^p-2aX^{\frac{p+1}{2}}+a^2X$ for some
  $a\in K$.
\item[(e)] $n=4$, $p=2$, and $f\sim_KX^4+(1+a)X^2+aX$ for some $a\in K$.
  \end{itemize}
  In the cases (b), (d), (e), and (a) for odd $n$, the following holds: If
  $K(w)$ is an intermediate field of $K(x,y)/K(t)$ with $[K(x,y):K(w)]=2$,
  then $K(w)$ is conjugate to $K(x)$.

  In case (a) suppose that $f(X)=D_n(X,a)$ and $K(w)$ is not conjugate to
  $K(x)$. Furthermore, suppose that $t=g(w)$ for a polynomial $g(X)\in K[X]$.
  Then $g(X)=-D_n(b(\xi+1/\xi)X+c,a)$ for $b,c\in K$ and $\xi$ a primitive
  $2n$-th root of unity.
\end{Proposition}
\begin{proof}
  Let $\hat K$ be the algebraic closure of $K$ in $K(x,y)$. Then
  $K(x)\subseteq \hat K(x)\subseteq K(x,y)$, so either $\hat K=K$ or
  $K(x,y)=\hat K(x)$.

  We start looking at the latter case. Here $\hat K(x)/\hat K(t)$ is a Galois
  extension with group $C$ which is a subgroup of $G=\Gal(\hat K(x)/K(t))$ of
  order $n$. Note that $C$ is either cyclic or dihedral. Let $\sigma\in C$, so
  $x^\sigma=\frac{ax+b}{cx+d}$ with $a,b,c,d\in\hat K$. From
  $f(\frac{ax+b}{cx+d})=f(x^\sigma)=f(x)^\sigma=t^\sigma=t=f(x)$ we obtain
  that $\frac{ax+b}{cx+d}$ is a polynomial, so $x^\sigma=ax+b$.

  Suppose that $p$ does not divide $n$. Then we may assume that the
  coefficient of $X^{n-1}$ of $f$ vanishes. From $f(ax+b)=f(x)$ we obtain
  $b=0$. Thus $C$ is isomorphic to a subgroup of $\hat K^\times$, in
  particular $C$ is cyclic and generated by $\sigma$ with $x^\sigma=\zeta x$
  with $\zeta$ a primitive $n$th root of unity. From $f(x)=f(\zeta x)$ we see
  that, up to a constant factor, $f(X)=X^n$. This is case (a) with $a=0$.

  From now on it is more convenient to work over an algebraic closure $\bar K$
  of $K$. As $\bar K(t)\cap K(x,y)=\hat K(t)$ (see e.g.\
  \cite[Prop.~1.11(c)]{Turnwald:Zakopane}), we obtain that $\Gal(\bar
  K(x)/\bar K(t))=C$.

  Now suppose that $p$ divides $n=\abs{C}$, but $p\ge3$. First assume that $C$
  is cyclic. From Lemma \ref{L:pglcyc} we get $p=n$. Let $\rho$ be a generator
  of $C$. Lemma \ref{L:pglcyc} shows the following: There is $x'\in\bar K(x)$
  with $\bar K(x)=\bar K(x')$, such that $x'^\rho=x'+1$. So $t'=x'^p-x'$ is
  fixed under $C$. We obtain $t'\in\bar K(t)$, because $\bar K(t)$ is the
  fixed field of $C$. From $p=[\bar K(x'):\bar K(t')]$ we obtain $\bar
  K(t')=\bar K(t)$. So there are rational functions $L,R\in\bar K(X)$ of
  degree $1$ with $x'=R(x)$ and $t=L(t')$. Then
  $f(x)=t=L(t')=L(x'^p-x')=L(r(x)^p-R(x))$, so $f=L\circ(X^p-X)\circ R$. By
  Lemma \ref{L:LRlr} we may assume that $L$ and $R$ are polynomials over $\bar
  K$. Then $f(X)=\alpha(X^p-aX)+\beta$ with $\alpha,\beta,a\in K$. From that
  we get case (b).

  Next assume that $C$ is dihedral of order $n$. As $p\ge3$, we get that $p$
  divides $n/2$. We apply Lemma \ref{L:pgldih} now. This yields $n=2p$, and
  there is $x'$ with $\bar K(x')=\bar K(x)$ such that $\bar K(t)$ is the fixed
  field of the automorphisms $x'\mapsto -x'$ and $x'\mapsto x'+1$. Obviously
  $t'=(x'^p-x')^2$ is fixed under these automorphisms, and as $[\bar
  K(x'):\bar K(t')]=2p$, we obtain $\bar K(t)=\bar K(t')$. The claim follows
  similarly as above.

  Now assume that $p=2$ divides $n$. Applying Lemmata \ref{L:pglcyc} and
  \ref{L:pgldih}, we get that $C$ is the Klein $4$ group. We see that
  $t'=x'(x'+1)(x'+b)(x'+b+1)$ is fixed under the automorphisms sending $x'$ to
  $x'+1$ and to $x'+b$. So $t'=h(x')$ with $h(X)=X^4+(1+b+b^2)X^2+(b+b^2)X$.
  Next we show that $b^2+b\in K$. A suitable substitution $\gamma f(\alpha
  X+\beta)+\delta$ should give $f(X)\in K[X]$.  We obtain $\gamma f(\alpha
  X+\beta)+\delta=\gamma(f(\alpha X)+f(\beta))+\delta\in K[X]$. Looking at the
  coefficients of $X^2$ and $X$ yields $\alpha\in K$, so $\alpha=1$ without
  loss of generality. Looking at $X^4$ gives $\gamma\in K$, so $\gamma=1$
  without loss. Finally the coefficient of $X$ yields the claim. Thus
  $f(X)=X^4+(1+b+b^2)X^2+(b+b^2)X\in K[X]$ and $\hat K=K(b)$, which gives case
  (e). In this case assume that $w$ is as in the proposition. Let $\tau_x$ and
  $\tau_w$ be the involutions of the dihedral group $G$ of order $8$ which fix
  $x$ and $w$, respectively. From $K(x,y)=K(x,b)=K(w,b)$ we obtain that
  $\tau_x,\tau_w\not\in C$. This shows that $\tau_x$ and $\tau_w$ are
  conjugate in $G$, so $K(w)$ is conjugate to $K(x)$.

  It remains to study the case $K=\hat K$, so $\hat K(x,y)/\hat K(t)$ is
  Galois with group $G$. By the Diophantine trick we obtain a rational
  parametrization of the quadric $F(X,Y)=0$ over $\bar K$ (actually, a
  suitable quadratic extension over which $F(X,Y)=0$ has a rational point
  suffices). In terms of fields that means $\bar K(z)=\bar K(x,y)$ for some
  element $z$.

  We apply Lemma \ref{L:pgldih}. Up to replacing $x$ and $t$ by $x'$ and $t'$
  as above, we get the following possibilities:

  (a) $p$ does not divide $n$, $x$ is fixed under the automorphism sending $z$
  to $1/z$, and $t$ is fixed under this automorphism and the one sending $z$
  to $z/\zeta$. So we may choose $t=z^n+1/z^n$, $x=z+1/z$. But then
  $t=D_n(x,1)$.  There are linear polynomials $L,R\in\bar K[X]$ with $L\circ
  D_n(X,1)\circ R=f\in K[X]$, so we get case (a) of the proposition by
  \cite[Lemma 1.9]{Turnwald:Schur}. For the remaining claims concerning this
  case, we may assume that $f(X)=D_n(X,a)$. Again set $t=f(x)$, and now choose
  $z$ with $z+a/z=x$. Then $t=D_n(x,a)=D_n(z+a/z,a)=z^n+(a/z)^n$. The normal
  closure $K(x,y)=K(x,w)$ of $K(x)/K(t)$ is contained in $K(\zeta,z)$. The
  elements $x'=\zeta x+\frac{a}{\zeta x}$ and $x''=\frac{x}{\zeta}+\frac{\zeta
    a}{x}$ are conjugates of $x$, so $x,x',x''\in K(x,y)$. From
  $x'+x''=(\zeta+1/\zeta)(x+a/x)$ we obtain $\zeta+1/\zeta\in K(x,y)$.
  However, we are in the case that $K$ is algebraically closed in $K(x,y)$, so
  $\zeta+1/\zeta\in K$.

  Suppose that $K(w)$ is not conjugate to $K(x)$. As extending the
  coefficients does not change Galois groups, this is equivalent to $\bar
  K(x)$ not being conjugate to $\bar K(w)$ in $\bar K(x,y)=\bar K(z)$. Note
  that $x$ is fixed under the involution $z\mapsto a/z$. The other involutions
  in $\Gal(\bar K(z)/\bar K(t))$ have the form $z\mapsto a\beta/z$, where
  $\beta$ is an $n$th root of unity, or $z\mapsto -z$. The latter involution
  cannot fix $w$, because the fixed field would be $\bar K(z^2)$, however,
  $z^n+(a/z)^n$ cannot be written as a polynomial in $z^2$. Thus suppose that
  $z\mapsto a\beta/z$ fixes $w$. If $\beta^{n/2}=1$, then an easy calculation
  shows that $\begin{bmatrix}0&a\\1&0\end{bmatrix}$ and
  $\begin{bmatrix}0&\beta a\\1&0\end{bmatrix}$ are conjugate in $\Gal(\bar
  K(z)/\bar K(t))$, contrary to $\bar K(x)$ and $\bar K(w)$ not being
  conjugate. Thus $\beta^{n/2}\ne1$, hence $\beta^{n/2}=-1$, because
  $\beta^n=1$. The element $w'=z+(\beta a)/z$ is fixed under the involution
  $z\mapsto a\beta/z$, so $\bar K(w')=\bar K(w)$. Furthermore,
  \[
  t=z^n+(a/z)^n=z^n+(\beta a/z)^n=D_n(z+(\beta a)/z,\beta a)=D_n(w',\beta a),
\]
so $g(X)=D_n(uX+v,\beta a)$ for some $u,v\in\bar K$. The condition that $g(X)$
has coefficients in $K$ shows that $\frac{v}{u}\in K$, see \cite[Lemma
1.9]{Turnwald:Schur}. Thus, upon replacing $X$ by $X-\frac{v}{u}$, we may
assume $v=0$. The transformation formula in Definition \ref{D:dickson} gives
$g(X)=D_n(uX,\beta
a)=\beta^{n/2}D_n(\frac{u}{\sqrt{\beta}}X,a)=-D_n(\frac{1}{\delta}X,a)$ with
$\delta\in\bar K$. As each conjugate of $w$ has degree $2$ over $K(x)$ we
obtain that $f(X)-g(Y)$ splits over $K$ in irreducible factors of degree $2$.
By Lemma \ref{L:dicksonfact} one of the factors of
$f(X)-g(Y)=D_n(X,a)+D_n(\frac{1}{\delta}Y,a)$ is
$X^2-\frac{1}{\delta}(\xi+1/\xi)XY+\frac{1}{\delta^2}Y^2-(\xi-1/\xi)^2a$. All
coefficients of this factor have to be in $K$, so there is $b_1\in K$ with
$\frac{1}{\delta}(\xi+\frac{1}{\xi})=b_1$. We obtain
$g(X)=-D_n(\frac{b_1}{\xi+1/\xi}X,a)=-D_n(b(\xi+1/\xi)X,a)$, where
$b=\frac{b_1}{(\xi+1/\xi)^2}\in K$. The claim follows.

(b) $n=p\ge3$. From a computation above we obtain $t=(z^p-z)^2$. We may assume
that $x$ is fixed under the automorphism sending $z$ to $-z$, so for instance
$x=z^2$. Let $h\in\bar K(X)$ with $h(x)=t$. That means
$h(z^2)=(z^p-z)^2=z^{2p}-2z^{p+1}+z^2$, hence $h(X)=X^p-2X^{\frac{p+1}{2}}+X$.
Lemma \ref{L:sim} yields the claim.

(c) The case $n=p=2$ does not arise, because we assumed $n\ge3$.

The conjugacy of $K(w)$ and $K(x)$ has been shown in the derivation of case
(e) above. In the cases (a) ($n$ odd), (b) and (d) it holds as well, because
$G$ is dihedral of order $2n$ with $n$ odd, so all involutions in $G$ are
conjugate.
\end{proof}
\section{Proof of the Theorems}
\subsection{Proof of Theorem \ref{T:main1} and \ref{T:main2}}
Suppose that $f(X)$ is not a polynomial in $X^p$, so not all coefficients of
$f$ are divisible by $p$. Let $q(X,Y)$ be an irreducible divisor of
$f(X)-g(Y)$ of degree at most $2$. Set $t=f(x)$, where $x$ is a transcendental
over $K$. Clearly both variables $X$ and $Y$ appear in $q(X,Y)$. In an
algebraic closure of $K(t)$ choose $y$ with $q(x,y)=0$. Note that $g(y)=t$.
The field $K(x)\cap K(y)$ lies between $K(x)$ and $K(t)$, so by L\"uroth's
Theorem, $K(x)\cap K(y)=K(u)$ for some $u$. Writing $t=\Phi(u)$ and $u=f_1(x)$
for rational functions $\Phi,f_1\in K(X)$, we have $f=\Phi\circ f_1$. By Lemma
\ref{L:LRlr}(a), we may replace $u$ by $u'$ with $K(u)=K(u')$, such that $t$
is a polynomial in $u$, and $u$ is a polynomial in $x$. Thus without loss of
generality we may assume that $\Phi$ and $f_1$ are polynomials. From that it
follows that $u$ is also a polynomial in $y$, so $g(X)=\Phi(g_1(X))$ for a
polynomial $g_1$ with $g_1(y)=u$. As $q$ is irreducible and
$f_1(x)-g_1(y)=u-u=0$, we get that $q(X,Y)$ divides $f_1(X)-g_1(Y)$. Thus, in
order to prove the theorems, we may assume that $f=f_1$ and $g=g_1$, so
$K(x)\cap K(y)=K(t)$.

First suppose that the polynomial $q(x,Y)$, considered in the variable $Y$, is
inseparable over $K(x)$. Then the characteristic of $K$ is $2$, and
$q(X,Y)=aX^2+bY^2+c$. This gives $x^2\in K(x)\cap K(y)=K(t)$, yielding case
(a) of Theorems \ref{T:main1} and \ref{T:main2}.

Thus we assume that $K(x,y)/K(x)$ is separable. By the assumption that $f(X)$
is not a polynomial in $X^p$ (this property is inherited by the new $f$), we
also obtain that $K(x)/K(t)$ is separable. Thus $K(x,y)/K(t)$ is separable.
From that one obtains the following: $K(x,y)/K(t)$ is Galois with group $G$,
and $G$ is generated by involutions $\tau_x$ and $\tau_y$, where $\tau_x$ and
$\tau_y$ fix $x$ and $y$, respectively. In particular, $G$ is a dihedral
group.

The case $\deg f=\deg g=2$ is trivial, thus assume $n=\deg f=\deg g\ge3$ from
now on.

The possibilities for $f$ are given in Proposition \ref{P:poldih}. In the
cases (b), (d), (e), and (a) for odd $n$, we obtain that $K(x)$ and $K(y)$ are
conjugate, yielding the case (a) of Theorem \ref{T:main1} and case (b) of
Theorem \ref{T:main2}.

Let us assume case (c) of Proposition \ref{P:poldih}. Here $G$ is a dihedral
group of order $4p$. If $\tau_x$ and $\tau_y$ are conjugate, then we obtain
case (a) of Theorem \ref{T:main1} and case (b)(iii) of Theorem \ref{T:main2}.
Thus suppose that $\tau_x$ and $\tau_y$ are not conjugate. By Lemma
\ref{L:dieder2} there is a conjugate $\tau_y'$ of $\tau_y$ such that $\tau_x$
and $\tau_y'$ generate a group of order $4$. Thus $K(x)$ and $K(y')$ have
degree $2$ over $K(x)\cap K(y')$. So there are $f_0,g_0,h\in K[X]$ with $f_0$
and $g_0$ of degree $2$ and $f=h\circ f_0$, $g=h\circ g_0$, giving case (a) of
Theorem \ref{T:main1}. Without loss of generality assume that
$f(X)=(X^p+aX+b)^2$, and $f_0(X)=X^2$. From $f(-X)=h((-X)^2)=h(X^2)=f(X)$ we
obtain $b=0$, so $f(X)=h(X^2)$ with $h(X)=X^p+2aX^{\frac{p+1}{2}}+a^2X$. This
yields case (d) of Theorem \ref{T:main2}.

Finally, assume the situation of Proposition \ref{P:poldih}, case (a) for even
$n$. If $K(x)$ and $K(y)$ are conjugate, then we obtain the case (a) of
Theorem \ref{T:main1} and case (b)(i) of Theorem \ref{T:main2}. If however
$K(x)$ and $K(y)$ are not conjugate, then Proposition \ref{P:poldih} yields
case (c) of Theorem \ref{T:main2}. In order to obtain case (b) of Theorem
\ref{T:main1} one applies Lemma \ref{L:dieder2} in order to show that $\tau_x$
and a conjugate of $\tau_y$ generate a dihedral $2$-group and argues as in the
previous paragraph.
\subsection{Proof of Theorem \ref{T:X^p}}
We have $f(X)=u(X)^p$ and $g(X)=v(X)^p$, where the coefficients of $u$ and $v$
are contained in a purely inseparable extension $L$ of $K$. (This includes the
case $K=L$.) In particular, $[L:K]$ is a power of $p$, so $q(X,Y)$ remains
irreducible over $L$ if $p>2$.

Suppose first that $p>2$, or that $q(X,Y)$ is irreducible over $L$ if $p=2$.
As each irreducible factor of $f(X)-g(Y)=u(X)^p-v(X)^p=(u(X)-v(Y))^p$ arises
at least $p$ times, we obtain that $q(X,Y)^p=q(X^p,Y^p)$ divides
$f(X)-g(Y)=f_0(X^p)-g_0(Y^p)$, and the claim follows in this case.

It remains to look at the case that $p=2$ and $q(X,Y)=q_1(X,Y)q_2(X,Y)$ is a
nontrivial factorization over $L$. If $q_1$ and $q_2$ do not differ by a
factor, then as above $q_1(X,Y)^2$ and $q_2(X,Y)^2$ divide $u(X)^2-v(Y)^2$, so
$q(X,Y)^2$ divides $u(X)^2-v(Y)^2$, and we conclude as above.

Thus $q(X,Y)=\delta(\alpha X+Y+\beta)^2$ for some $\alpha,\beta\in L$,
$\delta\in K$. Then $q(X,Y)=\delta(aX^2+Y^2+b)$ with $a,b\in K$ divides
$f_0(X^2)-g_0(Y^2)$, so $aX+Y+b$ divides $f_0(X)-g_0(Y)$, hence
$g_0(X)=f_0(aX+b)$, and the claim follows.
\begin{Remark}
  The method of the paper is easily extended to the study of degree $2$
  factors of polynomials of the form $a(X)b(Y)-c(X)d(Y)$, where $a,b,c,d$ are
  polynomials. For if $q(X,Y)$ is a quadratic factor, $x$ is a transcendental,
  and $y$ chosen with $q(x,y)=0$, then $a(x)/c(x)=d(y)/b(y)$, so setting
  $t=a(x)/c(x)=d(y)/b(y)$ and studying the field extension $K(x,y)/K(t)$
  requires only minor extensions of the arguments given in the paper.
\end{Remark}
%\bibliographystyle{myalpha}
%\bibliography{mybib}

\vskip.5cm
\noindent{\sc Poornaprajna Institute of Scientific Research, Davanhalli,
  Bangalore, India}\par\vspace{3mm}
\noindent{\sc Institut f\"ur Mathematik, Universit\"at W\"urzburg, Am Hubland,
  D-97074 W\"urzburg, Germany}\par
\noindent{\sl E-mail: }{\tt Peter.Mueller@mathematik.uni-wuerzburg.de}\par
\noindent{\sl URL: }{\tt
  www.mathematik.uni-wuerzburg.de/\~{}mueller}\par\vspace{3mm}
\noindent{\sc Statistics \& Mathematics Unit, Indian Statistical Institute,
  8th Mile Mysore Road, Bangalore -- 560 059}\par
\noindent{\sl E-mail: }{\tt sury@ns.isibang.ac.in}\par
\noindent{\sl URL: }{\tt www.isibang.ac.in/\~{}sury}
\end{document}